\documentclass[12pt]{article}
\usepackage[a4paper]{geometry}
\usepackage{amsthm}
\usepackage{amsmath}
\usepackage{amssymb}
\usepackage{amsfonts}
\usepackage{graphicx}
\usepackage{graphics}
\usepackage[bf,SL,BF]{subfigure}
\usepackage{srcltx}
\usepackage{multirow}
 \numberwithin{equation}{section}
\def\eref#1{(\ref{#1})}
\def\N{\mathbb{N}}

\def\R{\mathbb{R}}

\def\Rh{\mathcal{R}_h}
\def\Rht{\mathcal{R}_{h,t}}
\def\Rha{\mathcal{R}_{h,0}}

\def\A{\mathcal{A}}

\def\Q{\mathcal{Q}}
\def\<{\big\langle}
\def\>{\big\rangle}

\def\diiv{\operatorname{div}}

\def\Tr{\operatorname{Trace}}

\def\det{\operatorname{det}}

\def\esssup{{\operatorname{esssup}}}

\newtheorem{Lemma}{Lemma}[section]

\newtheorem{Theorem}{Theorem}[section]

\newtheorem{Condition}{Condition}[section]

\theoremstyle{remark}
\newtheorem{Remark}{Remark}[section]
\newtheorem{Example}{Example}[section]
\theoremstyle{definition}

\newtheorem{Algorithm}{Algorithm}[section]

\newtheorem{Assumption}{Assumption}[section]
\theoremstyle{definition}

\newcommand{\figbox}[1]{%
  \fbox{%
    \vbox to 1in{%
      \vfil
      \hbox to 2in{%
        \hfil
        #1%
        \hfil}%
      \vfil
    }
  }
}

\newcommand{\goodgap}{%
  \hspace{\subfigcapskip}}
\begin{document}
\title{Numerical Homogenization of the Acoustic Wave Equations with a Continuum of Scales.\protect\footnotetext{AMS 2000 {\it{Subject Classification}}. Primary 35L05, 
35B27; secondary 65M15, 86-08, 74Q15.}
\protect}
% Enter your title between
\date{\today}

\author{Houman Owhadi\footnote{California Institute of Technology
Applied \& Computational Mathematics, Control \& Dynamical systems,
MC 217-50 Pasadena , CA 91125, owhadi@caltech.edu}  and 
Lei Zhang\footnote{California Institute of Technology
Applied \& Computational Mathematics MC 217-50 Pasadena , CA 91125,
zhanglei@acm.caltech.edu}}

\maketitle
\begin{abstract}
In this paper, we consider numerical homogenization of acoustic wave 
equations with heterogeneous coefficients, namely, when the bulk modulus
and the density of the medium are only bounded. We show that under a 
Cordes type condition the second order derivatives of the solution with
respect to harmonic coordinates are $L^2$ (instead $H^{-1}$ with 
respect to Euclidean coordinates) and the solution itself is in 
$L^{\infty}(0,T,H^2(\Omega))$ (instead of $L^{\infty}(0,T,H^1(\Omega))$
with respect to Euclidean coordinates). Then, we propose an implicit 
time stepping method to solve the resulted linear system on coarse 
spatial scales, and present error estimates of the method. It follows 
that by pre-computing the associated harmonic coordinates, it is 
possible to numerically homogenize the wave equation without 
assumptions of scale separation or ergodicity.
\end{abstract}
\maketitle

\section{Introduction and main results}\label{sec1}

Let $\Omega$ be a bounded and convex domain of class $C^{2}$ of 
$\mathbb{R}^{2}$. Let $T>0$. Consider the following acoustic wave equation

\begin{equation}\label{waveeqn}
\begin{cases}
K^{-1}(x)\partial_t^2u=\diiv(\rho^{-1}(x)\nabla u(x,t))+g\quad\text{in}\quad\Omega\times(0,T).\\
u(x,t)=0\quad\text{for}\quad(x,t)\in\partial\Omega\times(0,T).\\
u(x,t)=u(x,0)\quad\text{for}\quad(x,t)\in\Omega\times\{ t=0\}.\\
\partial_t u(x,t)=u_t(x,0)\quad\text{for}\quad(x,t)\in\Omega\times\{ t=0\}.
\end{cases}
\end{equation}

Write $\Omega_T:=\Omega \times (0,T)$ and $a:=\rho^{-1}$. We assume $a$ 
is a uniformly elliptic $2\times 2$ symmetric matrix on $\bar{\Omega}$ whose 
entries are bounded and measurable. There exists $0<a_{\min}\leq a_{\max}$, 
such that $\forall\xi\in\R^2$, $|\xi|=1$, 
$a_{\min}\leq {^t\xi} a(x)\xi\leq a_{\max}$, $\forall x\in\Omega$. 
$K$ is a scalar such that $K_{\min}\leq K\leq K_{\max}$. 
$g\in L^{2}(\Omega_T)$.

Equation \eref{waveeqn} can be used to model wave propagation in heterogeneous 
media. It is important in many applications such as geophysics, seismology, and
electromagnetics \cite{MR1716789,  MR1442931, Wy98, VdMiKo05}. In geophysical 
and seismic prospecting, $K$ stands for the bulk modulus, $\rho$ the density 
and $u$ the unknown pressure. The velocity $c$ and acoustic impedance $\sigma$ 
are given by $c=\sqrt{K/\rho}$ and $\sigma=\sqrt{K\rho}$.

Wave propagation in heterogeneous media involves many different spatial
scales. Even with modern state-of-the-art supercomputers, a direct simulation 
of the highly heterogeneous media is often difficult if not impossible. That is 
why we want to use multiscale methods to solve \eref{waveeqn} on the coarse spatial 
scales. More precisely, we want to know how to transfer information from fine 
scales to coarse scales, and how to use the information obtained to solve the coarse 
scale problem with much fewer degrees of freedom. We often refer this procedure 
as numerical homogenization or numerical upscaling.  

The idea of using oscillating tests functions in relation to homogenization 
can be backtracked to the work of Murat and Tartar on homogenization and 
H-convergence, in particular we refer to \cite{MR1493039} and \cite{MR557520} 
(recall also that the framework of H-convergence is independent from 
ergodicity or scale separation assumptions). The implementation and practical 
application of oscillating test functions in finite element based numerical 
homogenization have been called multi-scale finite element methods (MsFEM) and 
have been studied by numerous authors \cite{AlBr05}, \cite{MR1286212}, 
\cite{MR1740386}, \cite{MR1232956}, \cite{MR1455261}. On 
the other hand, numerical schemes have been developed to solve the acoustic 
wave equation with discontinuous coefficients, for example in 
\cite{MR1442931} by nonconforming finite element method and in 
\cite{MR1716789} by domain decomposition. Recently, numerical homogenization 
or numerical upscaling methods such as \cite{VdMiKo05} are proposed for wave 
equation with heterogeneous coefficients. 

The finite element method in this paper is closer in spirit to the work of Hou 
and Wu \cite{MR1455261} and Allaire and Brizzi \cite{AlBr05}. It is based on a 
technique first introduced in \cite{OwZh05} for elliptic equations and 
extended in \cite{OwZh05b} to parabolic equations characterised by a continuum 
of scales in space and time. The main difference lies in the fact that, 
instead of solving a local cell problem to get a basis function as
in MsFEM (Multiscale finite element method) or to calculate effective media 
property as in upscaling method \cite{Far02}, we use a global change of 
coordinates. The global change of coordinates allows to avoid the so called 
cell resonance problem and obtain a scheme converging in situations where the 
medium has no separation between scales. This makes our method amenable to 
problems with strongly non-local medias, such as high conductivity channels. 

We use a composition rule to construct the finite element space. I. Babu\v{s}ka
et al. introduced the so called ``change of variable'' technique 
\cite{MR1286212} in the general setting of partition of unity method 
(PUM) with p-version of finite elements. Through this change of variable, the 
original problem is mapped into a new one which can be better approximated. 
Allaire and Brizzi \cite{AlBr05} introduced the composition rule in the 
multi-scale finite element formulation, and have observed that a multi-scale 
finite element method with higher order Lagrange polynomials has a higher 
accuracy. 

The main difference with parabolic equations \cite{OwZh05b} lies in the fact 
that with hyperbolic equations, energy is conserved and after homogenization 
there is no hope of recovering the energy (or information) lying in the 
highest frequencies. However when the medium is highly heterogeneous the 
eigenfunctions associated to the highest frequencies are localized, thus 
energy is mainly transported by the lowest frequencies. That is why, when 
one is only interested in the large scale transport of energy it is natural 
to approximate the solutions of \eref{waveeqn} by the solutions of an 
homogenized operator. For localization of waves in heterogeneous media, we 
refer to \cite{An58,MR1480976,Sh95}.

This paper is organized as follows. In the next section, we present the 
formulation of the mathematical problem and numerical methods, and also show 
main results. In Section 3, we will give the detailed proof and explanation 
of the results in Section 2. In Section 4, we present several numerical 
examples and conclusions.
  
\section{Main Results}\label{c4s1}
In general, the approximation power of finite element method is subject 
to the best approximation for an exact solution with respect to the finite 
element space. Therefore, we require smoothness of the solution to prove 
convergence theorems. That is one of the reasons why standard methods are not
applicable for problems with heterogeneous media. For example, in 
\eref{waveeqn}, we only have $u\in L^{\infty}(0,T,H^1(\Omega))$, and we 
can not gain anything if we approximate the solution with usual $C^0$ or $C^1$
finite element basis. However, as in \cite{OwZh05}, we can find harmonic 
coordinates which the solution of the wave equation is smoothly dependent on,
which is the so-called compensation phenomena.

\subsection{Compensation Phenomena}
We will focus on space dimension $n=2$. The extension to higher dimension is 
straightforward conditioned on the stability of $\sigma$. Let 
$F:=(F_{1},F_{2})$ be the harmonic coordinates satisfying 
\begin{equation}\label{dgdgfsghsf62c4}
\begin{cases}
\diiv a \nabla F=0 \quad \text{in}\quad \Omega,\\
F(x)=x \quad \text{on}\quad \partial \Omega.
\end{cases}
\end{equation}

Let $\sigma:={^t\nabla F}a\nabla F$ and 
\begin{equation}
\mu_{\sigma}:=\esssup_{x\in\Omega}\Big(\frac{\lambda_{\max}\big(\sigma(x)\big)}
{\lambda_{\min}\big(\sigma(x)\big)}\Big).
\end{equation}

\begin{Condition}\label{Cordes}
$\sigma$ satisfies Cordes type condition if: $\mu_\sigma<\infty$ and 
$(\Tr[\sigma ])^{-1}\in L^\infty(\Omega)$. 
\end{Condition}

\begin{Remark}\label{rmk21}
If $F$ is a quasiregular mapping, i.e., the dilation quotient (the ratio of 
maximal to minimal singular values of the Jacobi matrix) is bounded, then
Cordes type condition \ref{Cordes} is satisfied \cite{MR2001070}. A invertible
quasiregular mapping is called quasiconformal. In \cite{MR2001070} and 
references therein, invertibility of $F$ is proved for 
$a\in L^{\infty}(\Omega)$. Some sufficient conditions for $F$ being 
quasiconformal were also given, for example, $\det(a)$ is locally H\"older 
continuous. Unfortunately, a counterexample with checkerboard structure was 
proposed, and it can be shown that $\mu_{\sigma}$ is unbounded at the 
intersecting point, which is known in mechanics as stress concentration.
However, we will show that as a solution technique, the numerical methods 
proposed in this paper also works for the cases with stress concentration.
\end{Remark}

Let $L^2\big(0,T;H^1_0(\Omega)\big)$ be the Sobolev space associated to the norm
\begin{equation}
\|v\|_{L^2(0,T;H^1_0(\Omega))}^2:=\int_0^T
\big\|v(.,t)\big\|_{H^1_0(\Omega)}^2\,dt
\end{equation}

Also, we define the norm of the space $L^{\infty}(0,T,H^2(\Omega))$ by 
\begin{equation}
\|v\|_{L^{\infty}(0,T,H^2(\Omega))}=\text{esssup}_{0\leq t \leq
T}\Big(\int_{\Omega} \sum_{i,j} \big(\partial_i
\partial_j v(x,t)\big)^2 \,dx\Big)^{\frac{1}{2}}.
\end{equation}

We require the right hand side $g$, initial value $u(x,0)$ and $u_t(x,0)$ to be
smooth enough, which is a reasonable assumption in many applications. For 
example, we can made the following assumptions,
\begin{Assumption} \label{assumption1}
Assume that the $g$ satisfies $\partial_t g\in L^2(\Omega_T)$, 
$g\in L^\infty(0,T,L^2(\Omega))$, initial data $u(x,0)$ and $\partial_t u(x,0)$
satisfy $\partial_t u(x,0)\in H^1(\Omega)$ and 
$\nabla a(x) \nabla u(x,0)\in L^2(\Omega)$ or equivalently
$\partial^2_t u(x,0)\in L^2(\Omega)$.
\end{Assumption}

We have the following compensation theorem, 
\begin{Theorem}\label{thm2}
Suppose that Cordes condition \ref{Cordes} and Assumption \ref{assumption1} 
hold, then $u\circ F^{-1}\in L^\infty(0,T,H^2(\Omega))$ and
\begin{equation}
\begin{split}
\|u\circ F^{-1}\|_{L^{\infty}(0,T,H^2(\Omega))}\leq&C\big(\|g\|_{L^{\infty}
(0,T,L^2(\Omega))}+\|\partial_tg\|_{L^2(\Omega_T)}+\|\partial_tu(x,0)\|_{H^1(
\Omega)}\\
&+\|\partial^2_tu(x,0)\|_{L^2(\Omega)}\big).
\end{split}
\end{equation}
The constant $C$ can be written as
\begin{equation}
C=C(n,\Omega,K_{\min},K_{\max},a_{\min},a_{\max})\mu_\sigma
\big\|(\Tr[\sigma])^{-1}\big\|_{L^\infty(\Omega)}.
\end{equation}
\end{Theorem}

\begin{Remark} \label{rmk15}
We have gained one more order of integrability in the harmonic coordinates 
since in general $u\in L^{\infty}(0,T,H^1(\Omega))$. The condition 
$g\in L^2(\Omega_T)$ is sufficient to obtain Theorem \ref{thm2} and the 
following theorems. For the sake of clarity we have preferred to restrict 
ourselves to $g\in L^\infty(0,T,L^2(\Omega))$.
\end{Remark}

\subsection{Numerical Homogenization in Space}\label{sub2}
Suppose we have a quasi-uniform mesh. Let $X^h$ be a finite dimensional 
subspace of $H^1_0(\Omega)\cap W^{1,\infty}(\Omega)$ with the following 
approximation properties: There exists a constant $C_X$ such that
\begin{itemize}
\item
Interpolation property, i.e., for all $f\in H^2(\Omega)\cap H^1_0(\Omega)$
\begin{equation}\label{approp}
\inf_{v\in X^h} \|f-v\|_{H^1_0(\Omega)}\leq C_X h \|f\|_{H^2(\Omega)}.
\end{equation}

\item
Inverse Sobolev inequality, i.e., for all $v\in X^h$,
\begin{equation}\label{appsrop3}
 \|v\|_{H^2(\Omega)}\leq
C_X h^{-1} \|v\|_{H^1_0(\Omega)},
\end{equation}
and
\begin{equation}\label{appssswerop3}
 \| v\|_{H^1_0(\Omega)}\leq
C_X h^{-1} \|v\|_{L^2(\Omega)}.
\end{equation}
\end{itemize}

These properties are known to be satisfied when $X^h$ is a $C^1$ finite 
element space. One possibility is to use weighted extended B-splines (WEB) 
method developed by K. H\"ollig in \cite{MR1928544,MR1860269}, these elements 
are in general $C^1$-continuous. They are obtained from tensor products of 
one dimensional B-spline elements. The homogeneous Dirichlet boundary 
condition is satisfied by multiplying the basis functions with a smooth weight 
function $\omega$ which satisfies $\omega=0$ at the boundary . 

Write the solution space $V^h$ as
\begin{equation}
V^h:=\big\{\varphi\circ F(x)\, :\,\varphi\in X^h\big\}.
\end{equation}

\begin{Remark}\label{rmk23}
We prove all the following theoretical results by using exact $F$, however, in 
the numerical implementations, we have to use discrete solution $F_d$ to 
approximate $F$. A complete justification of the numerical scheme requires us to prove
$\|\phi\circ F-\phi\circ F_d\|\to 0$ as we refine the mesh. In \cite{AlBr05}, Allaire and 
Brizzi proved the convergence with respect to the discrete map $F_d$ in the periodic case 
using asymptotic expansion, as well as some regularity assumptions requiring the mappings $F$ 
and $F_d$ smooth enough. However, in the general case, neither the tool of 
asymptotic expansion nor smoothness assumption is available, which makes the
complete justification very difficult. Some discussions and further suggestions for similar 
problems in the context of variational mesh generation can be found in \cite{Ga04}. Another 
problem is, although $F$ is guaranteed to be invertible, $F_h$ is not. Fortunately this can be 
relieved if $F$ is solved by piecewise linear finite element and the mesh only has 
non-obtuse-angled triangles \cite{Fl03}. In view of the above discussion, we need to compute 
$F_d$ at a fine mesh such that the error $\|F-F_d\|$ is very small to fully resolve the small 
scale structure of $F$. 
\end{Remark}

We use the following notation
\begin{equation}\label{gaz52c4}
a[v,w]:=\int_\Omega {^t\nabla v(x,t)}a(x)\nabla w(x,t)\,dx.
\end{equation}

For $v\in H^1_0(\Omega)$ write $\Rha v$ the Ritz-Galerkin projection of 
$v$ on $V^h$ with respect to the bilinear operator $a[\cdot,\cdot]$, i.e., the
unique element of $V^h$ such that for all $w\in V^h$,
\begin{equation}
a[w,v-\Rha v]=0.
\end{equation}

Define $Y_T^h$ the subspace of $L^2\big(0,T;H^1_0(\Omega)\big)$ as
\begin{equation}
Y_T^h:=\{v\in L^2\big(0,T;H^1_0(\Omega)\big):v(x,t)\in V^h,\forall t\in[0,T]\}.
\end{equation}

Write $u_h$ the solution in $Y_T^h$ of the following system of ordinary 
differential equations:
\begin{equation}\label{ghjbbfh52c4}
\begin{cases}
(K^{-1}\psi(x),\partial_t^2 u_h)_{L^2(\Omega)}+a[\psi(x),u_h]=
(\psi(x),g)_{L^2(\Omega)}
\quad\text{for all $t\in(0,T)$ and $\psi\in V^h$},\\
u_h(x,0)=\Rha u(x,0),\\
\partial_tu_h(x,0)=\Rha\partial_tu(x,0).
\end{cases}
\end{equation}

The following theorem shows the error estimate of the semidiscrete solution. 
We need more smoothness on the forcing term $g$ and the initial data than 
Assumption \ref{assumption1} to guarantee the $O(h)$ convergence of the 
scheme \eref{ghjbbfh52c4}. On the other hand, we can see that even if $g$ 
and all the initial data are smooth, with general conductivity matrix $a(x)$,
we can merely expect $u\in L^{\infty}(0,T,H^1(\Omega))$ instead of the 
improved regularity $L^{\infty}(0,T,H^2(\Omega))$ in the harmonic 
coordinates, and the convergence rates will deteriote for the conventional 
finite elements.  

\begin{Assumption}\label{assumption2} 
Assume that the forcing term $g$ satisfies
$\partial^2_t g\in L^2(\Omega_T)$, $\partial_t g\in L^\infty(0,T,L^2(\Omega))$,
initial value $u(x,0)$ and $\partial_t u(x,0)$ satisfy
$\partial^2_t u(x,0)\in H^1(\Omega)$ and
$\nabla a(x) \nabla \partial_t u(x,0)\in L^2(\Omega)$ or equivalently
$\partial^3_t u(x,0)\in L^2(\Omega)$.
\end{Assumption}

From now on we will always suppose without explicitly mentioning that 
Assumption \ref{assumption2} is satisfied in the discussion of numerical 
homogenization method.

\begin{Theorem}\label{thlehsgsm5}
Suppose that Cordes condition \ref{Cordes} and Assumption \ref{assumption2} 
hold, we have
 \begin{equation}
 \begin{split}
\big\|\partial_{t}(u-u_{h})(.,T)\big\|_{L^2(\Omega)}+
\big\|(u-u_{h})(.,T)\big\|_{H^1_0(\Omega)}\leq C
 h\big(\|\partial_t g\|_{L^{\infty}(0,T,L^2(\Omega))}+
\|\partial^2_t g\|_{L^2(\Omega_T)}\\
+\|\partial^2_t u(x,0)\|_{H^1(\Omega)}
+\| \partial^3_t u(x,0) \|_{L^2(\Omega)} \big).
 \end{split}
 \end{equation}
The constant $C$ depends on $C_X$, $n$, $\Omega$, $\mu_\sigma$, $K_{\min}$, 
$K_{\max}$, $a_{\min}$, $a_{\max}$, and
$\big\|(\Tr[\sigma])^{-1}\big\|_{L^\infty(\Omega)}$. 
\end{Theorem}

\subsection{Numerical Homogenization in Time and Space}\label{dits}
Let $M\in\N$, $(t_{n}=n\frac{T}{M})_{0\leq n\leq M}$ is a discretization of 
$[0,T]$. $(\varphi_{i})$ is a $C^1$ basis of $X^h$. Write trial space $Z_T^h$ 
the subspace of $Y_T^h$ such that
\begin{equation}
\begin{split}
Z_T^h=\{w\in Y_T^h:w(x,t)=&\sum_{i}c_i(t)\varphi_i\big(F(x)\big),\mbox{
$c_i(t)$ are linear on $(t_n,t_{n+1}]$ and}\\
&\mbox{continuous on $[0,T]$}\}
\end{split}
\end{equation}

Let test space $U_T^h$ be the subspace of $Y_T^h$ such that
\begin{equation}
U_T^h=\{\psi\in Y_T^h:\psi(x,t)=\sum_{i} d_i \varphi_i\big(F(x)\big),\mbox{
$d_i$ are constant on $[0,T]$.}\}.
\end{equation}

Write $v_h$ the solution in $Z_T^h$ of the following system of implicit weak
formulation: for $n\in\{0,\dots,M-1\}$ and $\psi\in U_T^h$,
\begin{equation}\label{timestepeqn1}
\begin{split}
(K^{-1}\psi,\partial_{t}v_{h})(t_{n+1})-(K^{-1}\psi,\partial_{t}v_{h})(t_{n})
=\int_{t_{n}}^{t_{n+1}}(K^{-1}\partial_{t}\psi,\partial_{t}v_{h})dt\\
-\int_{t_{n}}^{t_{n+1}}a[\psi,v_{h}]dt+\int_{t_{n}}^{t_{n+1}}(\psi,g)dt.
\end{split}
\end{equation}

In equation \eref{timestepeqn1},  $\partial_{t}v_{h}(t)$ stands for
$\lim_{\epsilon\downarrow 0}(v_h(t)-v_h(t-\epsilon))/\epsilon$. 
Once we know the values of $v_h$ and $\partial_tv_h$ at $t_n$, 
\eref{timestepeqn1} is a linear system for the the unknown coefficients of 
$\partial_tv_h(t_{n+1})$ in $V^h$. By continuity of $v_h$ in time, we can 
obtain $v_h(t_{n+1})$ by
\begin{equation}\label{timestepeqn2}
v_h(t_{n+1})=\partial_{t}v_{h}(t_{n+1})({t_{n+1}-t_n})+v_h(t_n).
\end{equation}

The following Theorem \ref{thm15} shows the stability of the implicit scheme 
\eref{timestepeqn1}:
\begin{Theorem}\label{thm15}
Suppose that Cordes condition \ref{Cordes} and Assumption \ref{assumption2} 
hold, we have
\begin{equation}
\begin{split}
\|\partial_{t}v_{h}(.,T)\|_{L^2(\Omega)}+&\|v_{h}(.,T)\|_{H^1_0(\Omega)}\leq 
C\big(\|g\|_{L^2(\Omega_T)}+\|\partial_t u(x,0)\|_{L^2(\Omega)}\\
&+\|u(x,0)\|_{H^1(\Omega)}\big)+Ch \big(\|\partial_t g\|_{L^{\infty}(0,T,L^2(
\Omega))}+\|\partial^2_t g\|_{L^2(\Omega_T)}\\
&+\|\partial^2_t u(x,0)\|_{H^1(\Omega)}+\|\partial^3_t u(x,0)\|_{L^2(\Omega)}
\big).
\end{split}
\end{equation}
The constant $C$ depends on $a_{\min}$, $a_{\max}$, $K_{\min}$,
$K_{\max}$, and $T$.
\end{Theorem}

The following Theorem \ref{thm16} gives us the error estimate for the scheme 
\eref{timestepeqn1}.
\begin{Theorem}\label{thm16}
Suppose that Cordes condition \ref{Cordes} and Assumption \ref{assumption2} 
hold, we have
\begin{equation}
\begin{split}
&\big\|(\partial_{t}u_{h}-\partial_{t}v_{h})(.,T)\|_{L^2(\Omega)}+\big
\|(u_{h}-v_{h})(.,T)\big\|_{H^1(\Omega)}\leq C \Delta t (1+h^{-1})\\
&\big(\|\partial_t g\|_{L^{\infty}(0,T,L^2(\Omega))}+\|\partial^2_tg\|_{L^2(
\Omega_T)}+\|\partial^2_t u(x,0)\|_{H^1(\Omega)}+\|\partial^3_t u(x,0)\|_{L^2(
\Omega)}\big).
\end{split}
\end{equation}
The constant $C$ depends on $C_X$, $T$, $a_{\min}$, 
$a_{\max}$, $K_{\min}$, $K_{\max}$, $\mu_{\sigma}$,
and $\big\|(\Tr[\sigma])^{-1}\big\|_{L^\infty(\Omega)}$.
\end{Theorem}

\section{Proofs}\label{c4s2}
The proofs are organized into three subsections corresponding to the three 
subsections of section \ref{c4s1}. 

\subsection{Compensation Phenomena: Proof of Theorem \ref{thm2}}
\begin{Lemma}\label{lem1c4}
We have
\begin{equation}
\begin{split}
\|\partial_t^2u\|^{2}_{L^2(\Omega)}(T)+a[\partial_t u](T)\leq &C(T,
\frac{K_{\max}}{K_{\min}},K_{\max})\Big(a[\partial_t u](0)\\
&+\|\partial_t^2 u(x,0)\|^{2}_{L^2(\Omega)}+\|\partial_tg\|^{2}_{L^2(\Omega_T)}
\Big).
\end{split}
\end{equation}
\end{Lemma}
\begin{proof}
In case $a$ is smooth, differentiating \eref{waveeqn} with respect to t, we 
have
\begin{equation}
K^{-1} \partial_t^3 u-\diiv a\nabla \partial_t u=\partial_t g.
\end{equation}
multiplying by $\partial_t^2 u$, and integrating over $\Omega$, we obtain that
\begin{equation}
\frac{1}{2}\frac{d}{dt}\|K^{-\frac{1}{2}}\partial_t^2 u\|^{2}_{L^2(\Omega)}+
\frac{1}{2}\frac{d}{dt}a[\partial_t u]=(\partial_t g,\partial_t^2 u)_{L^2(
\Omega)}.
\end{equation}

Integrating the latter equation with respect to $t$ and using Cauchy-Schwartz 
inequality we obtain that
\begin{equation}
\begin{split}
\|K^{-\frac{1}{2}}\partial_t^2 u\|^{2}_{L^2(\Omega)}(T)+a[\partial_t u](T)\leq&
\|K^{-\frac{1}{2}}\partial_t^2 u\|^{2}_{L^2(\Omega)}(0)+a[\partial_t u](0)\\
&+\|\partial_t g\|_{L^2(\Omega_T)}\|\partial_t^2 u\|_{L^2(\Omega_T)}.
\end{split}
\end{equation}

Consider the following differential inequality, suppose that $A$ is constant, 
$B(t)>0$ and nondecrease, $X(t)>0$ and $X(t)$ is continuous with respect to $t$,
\begin{equation}
X(t) \leq A + B(t)\big(\int_0^t X(s)\,ds\big)^{\frac{1}{2}}.
\end{equation}

Write $Y(t)=\sup_{s\in [0,t]}X(s)$, one has
\begin{equation}
X(t)\leq A+B(t)t^\frac{1}{2}\big(Y(t)\big)^\frac{1}{2}\leq A+
\frac{t(B(t))^2+Y(t)}{2}.
\end{equation}

Take the supremum of both sides over $t\in[0,T]$, we have
\begin{equation}\label{diffineq}
Y(T)\leq 2A+T\big(B(T)\big)^2.
\end{equation}

It follows that
\begin{equation}
\begin{split}
\|\partial_t^2 u\|^{2}_{L^2(\Omega)}(T)+a[\partial_t u](T)\leq&C(T,
\frac{K_{\max}}{K_{\min}},K_{\max})\Big(a[\partial_t u](0)\\
&+\|\partial_t^2 u\|^{2}_{L^2(\Omega)}(0)+
\|\partial_t g\|^{2}_{L^2(\Omega_T)}\Big).
\end{split}
\end{equation}

In the case where $a$ is nonsmooth we use Galerkin approximations of $u$ in
\eref{waveeqn} and then pass to limit. This technique is standard and we refer
to \cite[Section 7.3.2.c]{Evans97} for a reminder.
\end{proof}

\begin{Lemma}\label{lem2}
\begin{equation}
\begin{split}
\|\partial_t u\|^{2}_{L^2(\Omega)}(T)+a[u](T)\leq& C(T,\frac{K_{\max}}
{K_{\min}},K_{\max})\Big(a[u](0)\\
&+\|\partial_t u\|^{2}_{L^2(\Omega)}(0)+\|g\|^{2}_{L^2(\Omega_T)}\Big).
\end{split}
\end{equation}
\end{Lemma}
\begin{proof}
Multiplying \eref{waveeqn} by $\partial_{t}u$, and integrating over $\Omega$, 
we obtain that 
\begin{equation}
\frac{1}{2}\frac{d}{dt}\|K^{-\frac{1}{2}} \partial_t u\|^{2}_{L^2(\Omega)}+
\frac{1}{2}\frac{d}{dt}a[u]=(g,\partial_t u)_{L^2(\Omega)}.
\end{equation}
The remaining part of the proof is similar to the proof of Lemma \ref{lem1c4}.
\end{proof}

We now need a variation of Campanato's result \cite{CM5} on non-divergence 
form elliptic operators. For a symmetric matrix $M$, let us write
\begin{equation}
\nu_M:=\frac{\Tr(M)}{\Tr({^tM M})}.
\end{equation}

Consider the following Dirichlet problem:
\begin{equation}\label{dcaqssaslkwaq21}
L_M v=f
\end{equation}

with $L_M:=\sum_{i,j=1}^2 M_{ij}(x) \partial_i \partial_j$. The
following  Theorem \ref{hdgjhdgd7} is an adaptation of Theorem 1.2.1 
of \cite{MPG00}. They are proved in \cite{MPG00} under the assumption that 
$M$ is bounded and elliptic. It can be proved that the conditions 
$\mu_{M}<\infty$ and $\nu_M<\infty$ are sufficient for the validity of the 
theorem, we refer to \cite{OwZh05b} and \cite{OwZh05}
for that proof.
\begin{Theorem}\label{hdgjhdgd7}
Assume that $\mu_M<\infty$,  $\nu_M\in L^\infty(\Omega)$ and $\Omega$
is convex. If $f\in L^2(\Omega)$ the Dirichlet problem
\eref{dcaqssaslkwaq21} has a unique solution satisfying
\begin{equation}\label{hdhazdgc7}
\|v\|_{W^{2,2}(\Omega)}\leq C \mu_M\|\nu_M f\|_{L^2(\Omega)}.
\end{equation}
\end{Theorem}

\begin{Remark}
The theorem can be extended to dimension $n>2$ under the general Cordes 
condition \cite{MPG00}.
\end{Remark}

Let us now prove the compensation result in Theorem \ref{thm2}. Choose
\begin{equation}\label{ksjshsj7622c4}
M:=\frac{\sigma}{|\det(\nabla F)|^\frac{1}{2}}\circ F^{-1}.
\end{equation}

Recall that $\sigma:={^t\nabla F}a\nabla F$. \eref{ksjshsj7622c4} is well 
defined since $\mu_M=\mu_\sigma$ and
\begin{equation}\label{hdhxcszswc7c4}
\|\nu_M\|_{L^\infty(\Omega)}^2 \leq
\frac{C}{\lambda_{\min}(a)}
\big\|(\Tr[\sigma])^{-1}\big\|^2_{L^\infty(\Omega)}.
\end{equation}

Fix $t\in [0,T]$. Choose
\begin{equation}\label{ghjhssalzz52c4}
\begin{split}
f:=\frac{(K^{-1}\partial_t^2 u-g)}{|\det(\nabla F)|^\frac{1}{2}}\circ F^{-1}.
\end{split}
\end{equation}

By the change of variable $y=F(x)$, one obtains that
\begin{equation}\label{ghjhssadsfalzz52c4}
\begin{split}
\|f\|_{L^2(\Omega)}\leq 2K_{\min}^{-1}\|\partial_t^2 u\|_{L^2(\Omega)}+
2\|g\|_{L^2(\Omega)}.
\end{split}
\end{equation}

Using the notation $\tilde{K}(y):=K(F^{-1}(y))$, 
$\tilde{g}(y,t):=g(F^{-1}(y),t)$, and $\tilde{u}(y,t):=u(F^{-1}(y),t)$, it 
follows from Theorem \ref{hdgjhdgd7} that there exists a unique 
$v\in W^{2,2}(\Omega)$ such that
\begin{equation}\label{ghjhsxfazz52c4}
\begin{split}
\sum_{i,j}\big(\sigma(F^{-1}(y))\big)_{i,j}\partial_i\partial_j v(y,t)=
\tilde{K}^{-1}(y)\partial_t^2 \tilde{u}(y,t)- \tilde{g}(y,t),
\end{split}
\end{equation}
and
\begin{equation}\label{hdhsszdgc7c4}
\|v\|_{W^{2,2}(\Omega)}\leq
C\mu_M\|\nu_M\|_{L^\infty(\Omega)} 
\big(K_{\min}^{-1}\|\partial_t^2 u\|_{L^2(\Omega)}+\|g\|_{L^2(\Omega)}\big).
\end{equation}

By change of variable $y=F(x)$ and the identity $\diiv a\nabla F=0$ we deduce 
that \eref{ghjhsxfazz52c4} can be written as
\begin{equation}\label{ghjhsxfsaazsdz52jc4}
\begin{split}
\diiv\big(a\nabla(v\circ F)\big)=K^{-1}\partial_t^2 u-g.
\end{split}
\end{equation}

If $\partial_t^2 u \in L^2(\Omega)$ and $g(.,t)\in L^2(\Omega)$ we can use the 
uniqueness property for the solution of the following divergence form 
elliptic equation (with homogeneous Dirichlet boundary condition)
\begin{equation}\label{ghjhsxfsaazsdz52c4}
\begin{split}
\diiv\big(a\nabla u\big)=K^{-1}\partial_t^2 u-g.
\end{split}
\end{equation}
to obtain that $v\circ F=u$. Thus we have proven Theorem \ref{thm2}. 

\subsection{Numerical Homogenization in Space: Proof of Theorem 
\ref{thlehsgsm5}.}
In the following sections we will prove the convergence of semidiscrete and fully discrete numerical homogenization formulation \eref{ghjbbfh52c4} and 
\eref{timestepeqn1}.

We have the following lemmas which are the discrete analogs of Lemma 
\ref{lem1c4} and Lemma \ref{lem4},
\begin{Lemma}\label{lem3c4}
We have
\begin{equation}
\begin{split}
\| \partial_t^2 u_h\|^{2}_{L^2(\Omega)}(T)+a[\partial_t u_h](T)\leq&
C(T,\frac{K_{\max}}{K_{\min}},K_{\max}) \Big( a[\partial_t u_h](0)\\
&+\| \partial_t^2 u_h(x,0)\|^{2}_{L^2(\Omega)}+
\|\partial_t g\|^{2}_{L^2(\Omega_T)}\Big).
\end{split}
\end{equation}
\end{Lemma}

\begin{Lemma}\label{lem4}
\begin{equation}
\begin{split}
\| \partial_t u_{h}\|^{2}_{L^2(\Omega)}(T)+a[u_h](T)\leq& 
C(T,\frac{K_{\max}}{K_{\min}},K_{\max}) \Big(a[u_h](0)\\
&+\|\partial_t u_h\|^{2}_{L^2(\Omega)}(0)+\|g\|^{2}_{L^2(\Omega)_T}\Big).
\end{split}
\end{equation}
\end{Lemma}

Write $\Rh$ the projection operator mapping $L^{2}(0,T;H_{0}^{1}(\Omega))$ 
onto $Y_T^h$, such that for all $v\in Y_T^h$:
\begin{equation}\label{proj}
\A_T[v,u-\Rh u]=0
\end{equation}
let $\rho:=u-\mathcal{R}_{h}u$ and $\theta:=\mathcal{R}_{h}u-u_{h}$, where
$u_h$ is the solution of \eref{ghjbbfh52c4}.

For fixed $t\in [0,T]$ and $v\in H^1_0(\Omega)$, we write 
$\mathcal{R}_{h,t}v(.,t)$ the solution of:
\begin{equation}
\int_\Omega {^t\nabla \psi}a(x)\nabla(v-\mathcal{R}_{h,t}v(x,t))\,dx=0
\quad\mbox{for all $\psi\in V^h$}
\end{equation}

It is obvious that $\mathcal{R}_h u(.,t)=\mathcal{R}_{h,t}u(.,t)$. For example,
we can choose a series of test functions in \eref{proj} which is separable in 
space and time, $v(x,t)=T(t)X(x)$, $T(t)$ is smooth in $t$ and has 
$\delta(t)$ function as its weak limit.

We need the following lemma:
\begin{Lemma}\label{lemrho} 
For $v\in H^1_0(\Omega)$ we have
\begin{equation}\label{ghhjsdaswwszwfh52}
\big(a[v-\Rht v]\big)^\frac{1}{2} \leq C h a_{\max}^\frac{1}{2}
\mu_\sigma^{\frac{1}{4}}\|\tilde{v}\|_{W^{2,2}}
\end{equation}
\end{Lemma}
\begin{proof}
Using  the change of coordinates $y=F(x)$ we obtain that (write
$\tilde{v}:=v \circ F^{-1}$)
\begin{equation}\label{ghhjasabbhfh52}
a[v]=Q[\tilde{v}]
\end{equation}
with
\begin{equation}\label{ghhssabbhfh52}
\Q[w]:=\int_\Omega {^t\nabla w(y)}Q(y)\nabla w(y)\,dy
\end{equation}
and
\begin{equation}\label{ghhsssh52}
Q(y):=\frac{\sigma}{\det(\nabla F)}\circ F^{-1}.
\end{equation}

Using the definition of $\Rh v$ we derive that 
\begin{equation}\label{ghhjasszbbhfh52}
\Q[\tilde{v}-\Rh v\circ F^{-1}]=\inf_{\varphi\in X^h}\Q[\tilde{v}-\varphi].
\end{equation}

By interpolation property \eref{approp} it follows,
\begin{equation}\label{ghhjawjbbhfh52}
\Q[\tilde{v}-\Rh v\circ F^{-1}]\leq \lambda_{\max}(Q) C_{X}^2 h^2
\|\tilde{v}\|_{W^{2,2}_D (\Omega)}^2.
\end{equation}
where $\lambda_{\max}(Q)$ is the supremum of eigenvalues of $Q$ over $\Omega$.

It is easy to obtain that
\begin{equation}
\lambda_{\max}(Q)\leq C a_{\max}\mu_\sigma^{\frac{1}{2}}
\end{equation}
which finishes the proof. 
\end{proof}

We will use the Lemmas \ref{lemrho}, \ref{tssnedhm2}, \ref{hgfg676},
\ref{slskj823c4} and \ref{lemdfdfd} to obtain the approximation property of
the projection operator $\mathcal{R}_h$.

With the improved Assumption \ref{assumption2}, differentiate \eref{waveeqn} 
with respect to $t$, and follow the proof of Theorem \ref{thm2}, we have 
\begin{Lemma}\label{tssnedhm2}
$\partial_t (u\circ F^{-1})\in L^{\infty}(0,T,H^2(\Omega))$ and
\begin{equation}
\begin{split}
\|\partial_t(u\circ F^{-1})\|_{L^{\infty}(0,T,H^2(\Omega))}\leq&C\big(\|
\partial_t g\|_{L^{\infty}(0,T,L^2(\Omega))}+\|\partial^2_t g\|_{L^2(\Omega_T)}
\\
&+\|\partial^3_t u(x,0)\|_{L^2(\Omega)}+\|\partial^2_t u(x,0)\|_{H^1(\Omega)}
\big).
\end{split}
\end{equation}
The constant $C$ is the one given in Theorem \ref{thm2}.
\end{Lemma}

Apply Lemma \ref{lemrho} to $\partial_t u$, we have
\begin{Lemma}\label{hgfg676}
\begin{equation}\label{sjswhawdsdlswec4}
\begin{split}
\big(\A_T[\partial_t \rho]\big)^{\frac{1}{2}}\leq&Ch\big(\|\partial_t 
g\|_{L^{\infty}(0,T,L^2(\Omega))}+\|\partial^2_t g\|_{L^2(\Omega_T)}\\
+&\|\partial^3_t u(x,0)\|_{L^2(\Omega)}+
\|\partial^2_t u(x,0)\|_{H^1(\Omega)}\big).
\end{split}
\end{equation}
The constant $C$ depends on $C_X$, $n$, $\Omega$, $\mu_\sigma$, 
$a_{\min}$, $a_{\max}$, $K_{\min}$,$K_{\max}$, and
$\big\|(\Tr[\sigma])^{-1}\big\|_{L^\infty(\Omega)}$ 
\end{Lemma}

We have the following estimate for $\|\partial_t\rho\|$ using the so-called 
Aubin-Nitsche trick \cite{Au67}. .
\begin{Lemma}\label{slskj823c4}
\begin{equation}\label{sjswhawdsdlswsswe}
\begin{split}
\|\partial_t \rho\|_{L^2(\Omega_T)} \leq &C h^2
\big(\|\partial_t g\|_{L^{\infty}(0,T,L^2(\Omega))}
+\|\partial^2_t g\|_{L^2(\Omega_T)}\\
 +& \|\partial^3_t u(x,0)\|_{L^2(\Omega)}+
\|\partial^2_t u(x,0)\|_{H^1(\Omega)}
\big).
\end{split}
\end{equation}
The constant $C$ in Lemma depends on $C_X$, $n$, $\Omega$, $\mu_\sigma$, 
$a_{\min}$, $a_{\max}$, $K_{\min}$, $K_{\max}$, and
$\big\|(\Tr[\sigma])^{-1}\big\|_{L^\infty(\Omega)}$ 
\end{Lemma}
\begin{proof}
We choose $v\in L^2(0,T,H^1_0(\Omega))$ to be the solution of the 
following linear
problem: for all $w\in L^2(0,T,H^1_0(\Omega))$
\begin{equation}\label{dhjdssh61}
A_T[w,v]=(w,\partial_t\rho)_{L^2(\Omega_T)}.
\end{equation}
Choosing $w=\partial_t\rho$ in equation \eref{dhjdssh61}  we deduce that
\begin{equation}
\|\partial_t\rho\|_{L^2(\Omega_T))}^2 =\A_T[\partial_t\rho,v-\Rh v].
\end{equation}
Using Cauchy Schwartz inequality we deduce that
\begin{equation}\label{sjhdlssswe}
\|\partial_t\rho\|_{L^2(\Omega_T)}^2\leq
\big(\A_T[\partial_t\rho]\big)^\frac{1}{2}
\big(\A_T[v-\Rh v]\big)^\frac{1}{2}.
\end{equation}
Since $\partial_t\rho(\cdot,t)\in L^2(\Omega)$, applying Theorem 
\ref{hdgjhdgd7} for $t\in[0,T]$ then integrate over $t$, we obtain that
\begin{equation}\label{ksjsjh61}
\|\hat{v}\|_{L^2(0,T,W^{2,2}(\Omega))}\leq C
\|\partial\rho\|_{L^2(\Omega_T)}.
\end{equation}
Using Lemma \ref{lemrho} we obtain that
\begin{equation}\label{sjhddsdlswe}
\big(\A_T[v-\Rh v]\big)^\frac{1}{2}\leq C h\|\partial_t\rho\|_{L^2(\Omega_T)}.
\end{equation}
It follows that
\begin{equation}\label{sjswhawdsdlsmxwe}
\|\partial_t\rho\|_{L^2(\Omega_T)}\leq C h 
\big(\A_T[\partial_t\rho]\big)^\frac{1}{2}.
\end{equation}
We deduce the lemma by applying  Lemma \ref{slskj823c4}  to
bound $A_T[\partial_t\rho]$.
\end{proof}

We have the following estimates for initial data,
\begin{Lemma}\label{lemdfdfd}
\begin{equation}
\begin{split}
\|\Rha\partial_t u(x,0)-\partial_t u(x,0)\|_{L^2(\Omega)}\leq Ch^2\big(\|\partial_tg(x,0)
\|_{L^2(\Omega)}+\|\partial^2_t u(x,0)\|_{H^1(\Omega)}+\|\partial^3_t u(x,0)
\|_{L^2(\Omega)}\big)\\
\|\Rha u(x,0)-u(x,0)\|_{H^1_0(\Omega)}\leq Ch\big(\|\partial_t g(x,0)\|_{L^2(
\Omega)}+\|\partial^2_t u(x,0)\|_{H^1(\Omega)}+\|\partial^3_t u(x,0)\|_{L^2(
\Omega)}\big)
\end{split}
\end{equation}
\end{Lemma}
\begin{proof}
We can estimate $\|\partial_t \rho\|_{L^2(\Omega)}$ using the duality argument
similar to Lemma \ref{slskj823c4} and derive the second inequality by Lemma 
\ref{lemrho}.
\end{proof}

\begin{Lemma}\label{lem5}
we have
\begin{equation}
\begin{split}
\|\partial_{t}(u-u_{h})\|^{2}_{L^2(\Omega)}(T)+a[u-u_{h}](T)\leq C(K_{\min},
K_{\max},T)\Big(\|\partial_{t}(u-u_{h})\|^{2}_{L^2(\Omega)}(0)\\
+a[u-u_{h}](0)+\|\partial_{t}\rho\|_{L^{2}(\Omega_T)}\|\partial_t^2(u-u_{h})
\|_{L^{2}(\Omega_T)}+\A_T[\partial_{t}\rho]\Big).
\end{split}
\end{equation}
\end{Lemma}
\begin{proof}
For $\psi\in L^2(0,T,H_0^1(\Omega))$, we have 
\begin{equation} 
(K^{-1}\psi,\partial_t^2(u-u_{h}))+a[\psi,u-u_{h}]=0.
\end{equation}

Let $\psi=\partial_t\theta=\partial_t(u-u_{h})-\partial_t\rho$, it follows
\begin{equation}
\frac{1}{2}\frac{d}{dt}\|K^{-\frac{1}{2}}\partial_{t}(u-u_{h})\|^{2}_{
L^2(\Omega)}+\frac{1}{2}\frac{d}{dt}a[u-u_{h}]=(K^{-1}\partial_{t}\rho,
\partial_t^2(u-u_{h}))+a[\partial_{t}\rho,u-u_{h}].
\end{equation}

Integrate with respect to $t$, using Cauchy-Schwartz inequality, we have
\begin{equation}
\begin{split}
&\frac{1}{2}\|K^{-\frac{1}{2}}\partial_{t}(u-u_{h})\|_{L^2(\Omega)}^2(T)-
\frac{1}{2}\|K^{-\frac{1}{2}}\partial_{t}(u-u_{h})\|_{L^2(\Omega)}^2(0)+
\frac{1}{2}a[u-u_{h}](T)\\
&-\frac{1}{2}a[u-u_{h}](0)\leq\int_{0}^{T}K^{-1}_{\min}\|\partial_{t}\rho\|_{
L^2(\Omega)}\|\partial_t^2(u-u_{h})\|_{L^2(\Omega)}dt+\big(\A_T[\partial_{t}
\rho]\A_T[u-u_{h}]\big)^\frac{1}{2}.
\end{split}
\end{equation}
The remaining part of the proof is similar to the proof of Lemma \ref{lem1c4}.
\end{proof}

Theorem \ref{thlehsgsm5} is a straightforward combination of Lemma
\ref{lem1c4}, Lemma \ref{lem3c4}, Lemma \ref{hgfg676}, Lemma \ref{slskj823c4}, 
Lemma \ref{lemdfdfd}, and Lemma \ref{lem5}.

\subsection{Numerical Homogenization in Space and Time: Proof of Theorems 
\ref{thm15} and \ref{thm16}}
\paragraph*{Stability}
Choose $\psi \in U_T^h$ in equation \eref{timestepeqn1} such that 
$\psi(x,t)=\partial_{t}v_{h}(x,t)$ for $t\in(t_n,t_{n+1}]$. We obtain that
\begin{equation}
\begin{split}
\|K^{-\frac{1}{2}}\partial_{t}v_{h}\|_{L^2(\Omega)}^{2}(t_{n+1})-(K^{-1}
\partial_{t}v_{h}(t_{n+1}),\partial_{t}v_{h}(t_{n}))_{L^2(\Omega)}=-
\int_{t_{n}}^{t_{n+1}}a[\partial_{t}v_{h},v_{h}]dt\\
+\int_{t_{n}}^{t_{n+1}}(\partial_{t}v_{h},g)_{L^2(\Omega)}dt.
\end{split}
\end{equation}

Observing that
\begin{equation}
\int_{t_{n}}^{t_{n+1}}a[\partial_{t}v_{h},v_{h}]dt=\frac{1}{2}a[v_{h}](t_{n+1})
-\frac{1}{2}a[v_{h}](t_{n}).
\end{equation}
 
using Cauchy-Schwartz inequality it follows,
\begin{equation}
\begin{split}
\|K^{-\frac{1}{2}}\partial_{t}v_{h}\|^{2}(t_{n+1})+a[v_{h}](t_{n+1})\leq
\|K^{-\frac{1}{2}}\partial_{t}v_{h}\|^{2}(t_{n})+a[v_{h}](t_{n})\\
+2\int_{t_{n}}^{t_{n+1}}(\partial_{t}v_{h},g)_{L^2(\Omega)}(t)\,dt.
\end{split}
\end{equation}

Summing over n from $0$ to $M-1$, we have,
\begin{equation}
\|K^{-\frac{1}{2}}\partial_{t}v_{h}\|^{2}(T)+a[v_{h}](T)\leq\|K^{-\frac{1}{2}}
\partial_{t}v_{h}\|^{2}(0)+a[v_{h}](0)+
2\int_{0}^{T}(\partial_{t}v_{h},g)_{L^2(\Omega)}\,dt.
\end{equation}

We conclude the proof of Theorem \ref{thm15} using the inequality
\eref{diffineq} in the proof of Lemma \ref{lem1c4}.

\paragraph*{$H^1$ Error Estimate}
We derive from equations \eref{timestepeqn1} and \eref{ghjbbfh52c4} that
\begin{equation}\label{kjsjshj82}
\begin{split}
(K^{-1}\psi,\partial_{t}u_{h}-\partial_{t}v_{h})(t_{n+1})-(K^{-1}\psi,
\partial_{t}u_{h}-\partial_{t}v_{h})(t_{n})\\
-\int_{t_{n}}^{t_{n+1}}(K^{-1}\partial_{t}\psi,\partial_{t}u_{h}-\partial_{t}
v_{h})dt+\int_{t_{n}}^{t_{n+1}}a[\psi,u_{h}-v_{h}]dt=0.
\end{split}
\end{equation}

Let $\psi=\partial_t\hat{u}_h-\partial_tv_h$ where $\hat{u}_h$ is the 
linear interpolation of $u_h$ over $Z_T^h$. Write $y_h=u_h-v_h$ and 
$w_h=\hat{u}_{h}-u_h$, it follows that
\begin{equation}\label{skjshs782}
\begin{split}
(K^{-1}\partial_{t}y_{h},\partial_{t}y_{h})(t_{n+1})+(K^{-1}\partial_{t}w_{h},
\partial_{t}y_{h})(t_{n+1})-(K^{-1}\partial_{t}y_{h},\partial_{t}y_{h})(t_{n})
\\
-(K^{-1}\partial_tw_h,\partial_ty_h)(t_n)+\int_{t_n}^{t_{n+1}}a[
\partial_ty_h,y_h]dt+\int_{t_n}^{t_{n+1}}a[\partial_tw_h,y_h]dt=0.
\end{split}
\end{equation}

Observing $\int_{t_n}^{t_{n+1}}\partial_{t}w_{h}(x,t)\,dt=0$ we
need the following lemma, which is a slight variation of the Hilbert-Bramble 
lemma,
\cite{BrHi70},
\begin{Lemma}\label{ineq}
If $\int^{t_{n+1}}_{t_n} u(s)\,ds=0$, then
\begin{equation}
u^2\leq\frac{1}{4}\Delta t\int^{t_{n+1}}_{t_n} u'(s)^2\,ds.
\end{equation}
\end{Lemma}

Since $\partial_t^2 w_h(x,t)=-\partial_t^2 u_h(x,t)$ in $(t_n,t_{n+1}]$, by 
Lemma \ref{ineq} we have
\begin{equation}\label{jshsg71}
\int_{\Omega}|\partial_t w_h(x,t)|^2\,dxdt\leq\frac{1}{4}\Delta t
\int_{t_n}^{t_{n+1}}\int_{\Omega}|\partial_t^2 u_h(x,t)|^2\,dx\,dt,
\end{equation} 
and
\begin{equation}\label{jshsg72}
\int_{t_n}^{t_{n+1}}\int_{\Omega}|\partial_t w_h(x,t)|^2\,dxdt\leq\frac{1}{4} 
\Delta t^2\int_{t_n}^{t_{n+1}}\int_{\Omega}|\partial_t^2 u_h(x,t)|^2\,dx\,dt.
\end{equation} 

Using the inverse Sobolev inequality \eref{appssswerop3} we obtain
from equation \eref{jshsg72} that
\begin{equation}\label{jshsgsshk72}
\int_{t_n}^{t_{n+1}}\int_{\Omega} a[\partial_t w_h]\,dx\,dt\leq C\frac{\Delta 
t^2}{h^2}\int_{t_n}^{t_{n+1}}\int_{\Omega} |\partial_t^2 u_h(x,t)|^2 \,dx\,dt.
\end{equation}

Summing \eref{skjshs782} over $n$, notice $y_h(0)=0$, $\partial_t y_h(0)=0$ we 
obtain that
\begin{equation}\label{kjsjsszhj82}
\begin{split}
(K^{-1}\partial_{t}y_{h},\partial_{t}y_{h})_{L^2(\Omega)}(T)+\frac{1}{2}a[y_{h}
(.,T)]=-\int_{0}^{T}a[\partial_{t}w_{h},y_{h}]dt-(K^{-1}(\partial_{t}w_{h},
\partial_{t}y_{h})_{L^2(\Omega)}(T).
\end{split}
\end{equation}

Theorem \ref{thm16} is a straightforward consequence of \eref{kjsjsszhj82}, 
the estimates \eref{jshsg71}, \eref{jshsgsshk72}, Lemma \ref{lem3c4} and Lemma 
\ref{lemdfdfd}.

\section{Numerical Experiments}\label{c4s3}
In this section, we will present the numerical algorithm and examples. 

We use web extended B-spline based finite element \cite{MR1928544} to span the 
space $X^h$ introduced in subsection \ref{sub2}. For all the numerical 
examples, we compute the solutions up to time $T=1$. The initial condition is 
$u(x,0)=0$ and $u_t(x,0)=0$. The boundary condition is $u(x,t)=0$, 
$x\in\partial\Omega$. For simplicity, the computational domain is the square 
$[-1,1]\times[-1,1]$ in dimension two.

We have a fine mesh and a coarse mesh characterized by different degrees of 
freedom (\emph{dof}). In general, the fine mesh is generated by hierarchical 
refinement of the coarse mesh: for each triangle of the coarse mesh, choose 
middle points of its 3 edges as new vertices, and divide the triangle into 4 new 
triangles. $a$ is defined as a piecewise constant function over each fine mesh
triangle, and is evaluated at the center of mass of the triangle.

\begin{Algorithm}[Algorithm for Numerical Homogenization]
\hfill\par
\begin{tabbing}
1. \=Compute F on fine mesh, the fine mesh solver for $F$ is \textit{Matlab} 
routine \textit{assempde}.\\
2. Construct multi-scale finite element basis $\psi=\varphi\circ F$, compute 
stiffness matrix $K$ and\\
\>mass matrix $M$.\\
3. March \eref{timestepeqn1} and \eref{timestepeqn2} in time with respect to 
the coarse \emph{dof}.\\
4. Repeat 3 if we have multiple right hand sides.
\end{tabbing}
\end{Algorithm}

In the implementation, $F$ is approximated by a piecewise linear finite element
solution. We mesh the square such that no triangle has an obtuse angle, 
therefore $F$ is an invertible piecewise linear mapping \cite{Fl03}. 
When we construct $\psi$, we simply take its piecewise linear interpolation
on the fine mesh.   

All the computations were done at a single Opteron Dual-Core 2600 cpu of a Sun Fire X4600 
server, and programmed in Matlab 7.3.

\begin{Example}
\label{exa:trig} Multiscale trigonometric coefficients
\end{Example}
The following example is extracted from \cite{MiYu06} as a problem without 
scale separation:
\begin{equation}
\begin{split}
a(x)=&\frac{1}{6}\Big(\frac{1.1+\sin(2\pi x/\epsilon_{1})}{1.1+\sin(2\pi
y/\epsilon_{1})}+\frac{1.1+\sin(2\pi y/\epsilon_{2})}{1.1+\cos(2\pi
x/\epsilon_{2})}+\frac{1.1+\cos(2\pi x/\epsilon_{3})}{1.1+\sin(2\pi
y/\epsilon_{3})}\\
&+\frac{1.1+\sin(2\pi y/\epsilon_{4})}{1.1+\cos(2\pi x/\epsilon_{4})}
+\frac{1.1+\cos(2\pi x/\epsilon_{5})}{1.1+\sin(2\pi y/\epsilon_{5})}+
\sin(4x^{2}y^{2})+1\Big)
\end{split}
\end{equation}
where $\epsilon_{1}=\frac{1}{5}$, $\epsilon_{2}=\frac{1}{13}$, 
$\epsilon_{3}=\frac{1}{17}$, $\epsilon_{4}=\frac{1}{31}$, 
$\epsilon_{5}=\frac{1}{65}$.
The conductivity $a$ is smooth, therefore it satisfies Cordes condition 
\ref{Cordes}.

First, we want to compare the performance of different numerical 
homogenization methods, 
\begin{itemize}
\item LFEM: A multi-scale finite element where $F$ is computed locally 
(instead of globally) on each triangle $K$ of the coarse mesh as the 
solution of a cell problem with boundary condition $F(x)=x$ on $\partial K$. 
This method has been implemented in order to understand the effect of the 
removal of global information in the structure of the metric induced by $F$.
\item FEM\_$\psi_{lin}$: The Galerkin scheme using the finite elements
$\psi_i=\varphi_{i}\circ F$, where $\varphi_i$ are the piecewise linear nodal 
basis elements.
\item FEM\_$\psi_{sp}$: The Galerkin scheme using the  finite element 
$\psi_i=\varphi_{i}\circ F$, where $\varphi_i$ are weighted cubic B-spline elements.
\end{itemize}

Suppose $u_f$ is the finite element solution of \eref{waveeqn} computed on
the fine mesh at time $T=1$, the fine mesh solver is \textit{Matlab} routine 
\textit{hyperbolic}, 
which uses linear finite element basis in space and adaptive ODE integrator 
in time. $v_h$ is the solution of \eref{timestepeqn1}. Numerical errors in
the norm $\|\cdot\|$ are computed by 
\begin{equation}
  error=\frac{\|v_h-u_f\|}{\|u_f\|}
\end{equation}
Numerical errors in $L^1$, $L^2$, $L^{\infty}$ and $H^1$ norm are computed.

\begin{table}
\begin{center}
\caption{Example \ref{exa:trig}, numerical errors of different methods, coarse \emph{dof} 49, fine \emph{dof} 261121, $g=1$} \label{ferrtindepp7c4}
\begin{tabular}{|c|c|c|c|c|}
\hline Method& $L^{1}$& $L^{\infty}$& $L^{2}$& $H^{1}$ \\
\hline
LFEM& 0.0440& 0.0982& 0.0534& 0.2054 \\
FEM\_$\psi_{lin}$& 0.0315& 0.0518& 0.0362& 0.1601\\
FEM\_$\psi_{sp}$& 0.0021& 0.0035& 0.0022& 0.0189\\
\hline
\end{tabular}
\end{center}
\end{table}

In Table \ref{ferrtindepp7c4} performances of different methods with coarse 
mesh \emph{dof} 49 are compared. We observe that the methods using global $F$ 
have better performance, and FEM\_$\psi_{sp}$ is much better than other 
methods. Note that the improvement of FEM\_$\psi_{lin}$ over LFEM is not
as significant as the elliptic case \cite{OwZh05}. 

From now on, all the results are computed by the method FEM\_$\psi_{sp}$.

Next, the impact of right hand side on accuracy will be investigated. We 
solve equation \eref{waveeqn} with a time independent source term $g=1$, a
slowly varying term $g=\sin(2.4x-1.8y+2\pi t)$, and a Gaussian source term 
given by
\begin{equation}\label{gaussian}
g(x,y)=\frac{1}{\sqrt{2\pi\sigma^2}}\exp\big(-\frac{x^2+(y-0.15)^2}{2\sigma^2}
\big)
\end{equation}
with $\sigma=0.05$. Notice that as $\sigma\rightarrow 0$, the source function 
will become singular in space. 

\begin{table}
\begin{center}
\caption{Example \ref{exa:trig}, numerical errors of FEM\_$\psi_{sp}$, with $g=1$, $dof_f$ is fine mesh \emph{dof}, $dof_c$ is coarse mesh \emph{dof}.}
\label{trigg1}
\begin{tabular}{|c|c|c|c|c|c|}
\hline $dof_f$&$dof_c$& $L^{1}$& $L^{\infty}$& $L^{2}$& $H^{1}$\\
\hline \multirow{3}{*}{65025} 
& 9& 0.0075& 0.0118& 0.0074& 0.0394 \\
& 49& 0.0023& 0.0037& 0.0023& 0.0194\\
& 225& 0.0009& 0.0023& 0.0010& 0.0117\\
\hline \multirow{3}{*}{261121} 
& 9& 0.0070& 0.0106& 0.0069& 0.0373 \\
& 49& 0.0021& 0.0035& 0.0022& 0.0188\\
& 225& 0.0009& 0.0025& 0.0010& 0.0117\\ \hline
\end{tabular}
\end{center}
\end{table}

\clearpage
\begin{table}
\begin{center}
\caption{Example \ref{exa:trig}, numerical errors of FEM\_$\psi_{sp}$, with
$g=\sin(2.4x-1.8y+2\pi t)$}
\label{triggsin}
\begin{tabular}{|c|c|c|c|c|c|}
\hline $dof_f$&$dof_c$& $L^{1}$& $L^{\infty}$& $L^{2}$& $H^{1}$\\
\hline \multirow{3}{*}{65025} 
& 9& 0.0400& 0.0390& 0.0360& 0.0869 \\
& 49& 0.0107& 0.0105& 0.0096& 0.0393\\
& 225& 0.0035& 0.0047& 0.0033& 0.0233\\
\hline \multirow{3}{*}{261121} 
& 9& 0.0399& 0.0373& 0.0359& 0.0866 \\
& 49& 0.0104& 0.0109& 0.0095& 0.0391\\
& 225& 0.0034& 0.0047& 0.0033& 0.0231\\ \hline
\end{tabular}
\end{center}
\end{table}

\begin{table}
\begin{center}
\caption{Example \ref{exa:trig}, numerical errors of FEM\_$\psi_{sp}$, with 
the Gaussian source $g$ in \eref{gaussian}}
\label{triggg}
\begin{tabular}{|c|c|c|c|c|c|}
\hline $dof_f$&$dof_c$& $L^{1}$& $L^{\infty}$& $L^{2}$& $H^{1}$\\
\hline \multirow{3}{*}{65025} 
& 9& 0.0581& 0.2270& 0.0704& 0.3484 \\
& 49& 0.0272& 0.1023& 0.0333& 0.2305\\
& 225& 0.0096& 0.0179& 0.0095& 0.0957\\
\hline \multirow{3}{*}{261121} 
& 9& 0.0574& 0.2199& 0.0688& 0.3436 \\
& 49& 0.0274& 0.976& 0.0332& 0.2254\\
& 225& 0.0097& 0.0212& 0.0101& 0.1005\\ \hline
\end{tabular}
\end{center}
\end{table}

Table \ref{trigg1} presents the errors for time independent term $g=1$. Table 
\ref{triggsin} presents the errors for slowly varying source
term $g=\sin(2.4x-1.8y+2\pi t)$. Table \ref{triggg} presents the errors for 
relatively singular Gaussian forcing term. It is clear that if the source term
is time independent and smooth in space, the method is more accurate, which
corresponds to the smoothness requirement of $g$ in Theorem \ref{thm16}. 
In all the examples, we have tried two fine meshes which have \emph{dof} 
$65025$ and $261121$ respectively, roughly $250\times 250$ and $500\times500$.
It can be seen that for fixed coarse \emph{dof}, the errors with respect to
different fine \emph{dof} are pretty close, which means fine \emph{dof} 
$65025$ is enough for this problem.

\begin{figure}[httb]
\begin{center}
\includegraphics[scale=0.3]{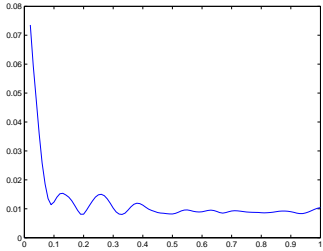}
\caption{Example \ref{exa:trig}, $L^1$ error with respect to time, for $g=1$, coarse \emph{dof} 9, fine \emph{dof} 65025}\label{terr}
\end{center}
\end{figure}

The Figure \ref{terr} shows the $L^1$ error evolution with respect to time, 
which is typical for other norms. The overshoot at the beginning corresponds
to the time discretization step. After several steps, the errors tend to be
stable.

\clearpage

\begin{figure}[httb]
\begin{center}
\includegraphics[scale=0.3]{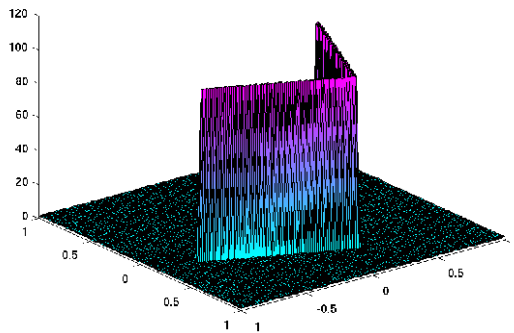}
\caption{Example \ref{exa:chanc4}, high conductivity channel medium}
\label{ap6c4}
\end{center}
\end{figure}

\begin{Example}
\label{exa:chanc4}Time independent high conductivity channel
\end{Example}
High conductivity channel is an interesting test problem in many petroleum 
applications because of its strong non-local effects. In this example $a$ is 
characterized by a narrow and long ranged high conductivity channel. We 
choose $a(x)=A\gg1$, if $x$ is in the channel, and $a(x)=1$, if $x$ is not in 
the channel. The media is illustrated in Figure \ref{ap6c4}. However, 
in this case, whether or not Cordes condition \ref{Cordes} is not clear.
We will go ahead testing the numerical performance of our method. 

Table \ref{ferrp6tindepspc4} shows numerical errors for $g=1$ with fixed coarse
\emph{dof} 49 and $A=10^1, 10^2, 10^3, 10^4$ respectively. From the table we 
can see that the errors grow with the aspect ratio increasing, but the growth 
is quite moderate and the numerical behavior of the method is stable. The 
errors for time dependent right hand side $g=\sin(2.4 x-1.8 y+2 \pi t)$ with 
$A=10^2 $ are also given in Table \ref{ferrtindeptdrhsp6}. 

\begin{table}
\begin{center}
\caption{Example \ref{exa:chanc4}, numerical errors with respect to different
aspect ratios, coarse \emph{dof} 49, fine \emph{dof} 261121.} \label{ferrp6tindepspc4}
\begin{tabular}{|c|c|c|c|c|}
\hline $A$& $L^{1}$& $L^{\infty}$& $L^{2}$& $H^{1}$\\
\hline
10& 0.0021& 0.0056& 0.0025& 0.0240\\
100& 0.0118& 0.0497& 0.0180& 0.0964\\
1000& 0.0181& 0.0931& 0.0316& 0.1308\\
10000& 0.0243& 0.1174& 0.0419& 0.1550\\
\hline
\end{tabular}
\end{center}
\end{table}

\begin{table}
\begin{center}
\caption{Example \ref{exa:chanc4}, numerical errors for $g=\sin(2.4x-1.8y+2\pi t)$}
\label{ferrtindeptdrhsp6}
\begin{tabular}{|c|c|c|c|c|c|}
\hline $dof_f$&$dof_c$& $L^{1}$& $L^{\infty}$& $L^{2}$& $H^{1}$\\
\hline \multirow{3}{*}{65025} 
& 9& 0.0750& 0.0777& 0.0729& 0.1528 \\
& 49& 0.0301& 0.0351& 0.00298& 0.0779\\
& 225& 0.0096& 0.0118& 0.0092& 0.0324\\
\hline \multirow{3}{*}{261121} 
& 9& 0.0752& 0.0779& 0.0731& 0.1533 \\
& 49& 0.0302& 0.0345& 0.0299& 0.0782\\
& 225& 0.0094& 0.0116& 0.0091& 0.0321\\ \hline
\end{tabular}
\end{center}
\end{table}

\begin{figure}[httb]
\begin{center}
\includegraphics[scale=0.3]{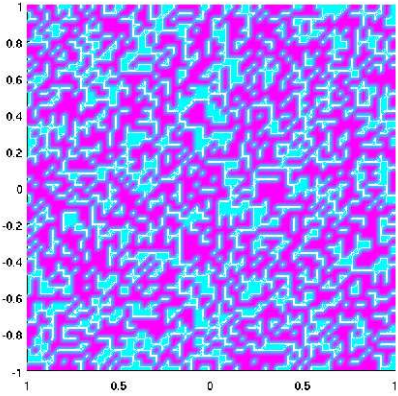}
\caption{Example \ref{exa:siteprcoc4}, site percolation medium} \label{ap7c4}
\end{center}
\end{figure}

\begin{Example}
\label{exa:siteprcoc4}Time independent site percolation
\end{Example}
In this example we consider the site percolating medium associated to Figure
\ref{ap7c4}. In this case, we subdivide the square into a $64\times 64$ checkerboard, the 
conductivity of each site is equal to $\gamma$ or $1/\gamma$ with probability $1/2$. We have chosen 
$\gamma=10$ in this example. In fact, this medium may not satisfy the Cordes condition 
\ref{Cordes} (also refer to Remark \ref{rmk21}). However, we will show that the method still 
works fine for this example.

\begin{figure}[httb]
  \begin{center}
    \subfigure[$u$.]
    {\includegraphics[width=0.35\textwidth,height= 0.3\textwidth]{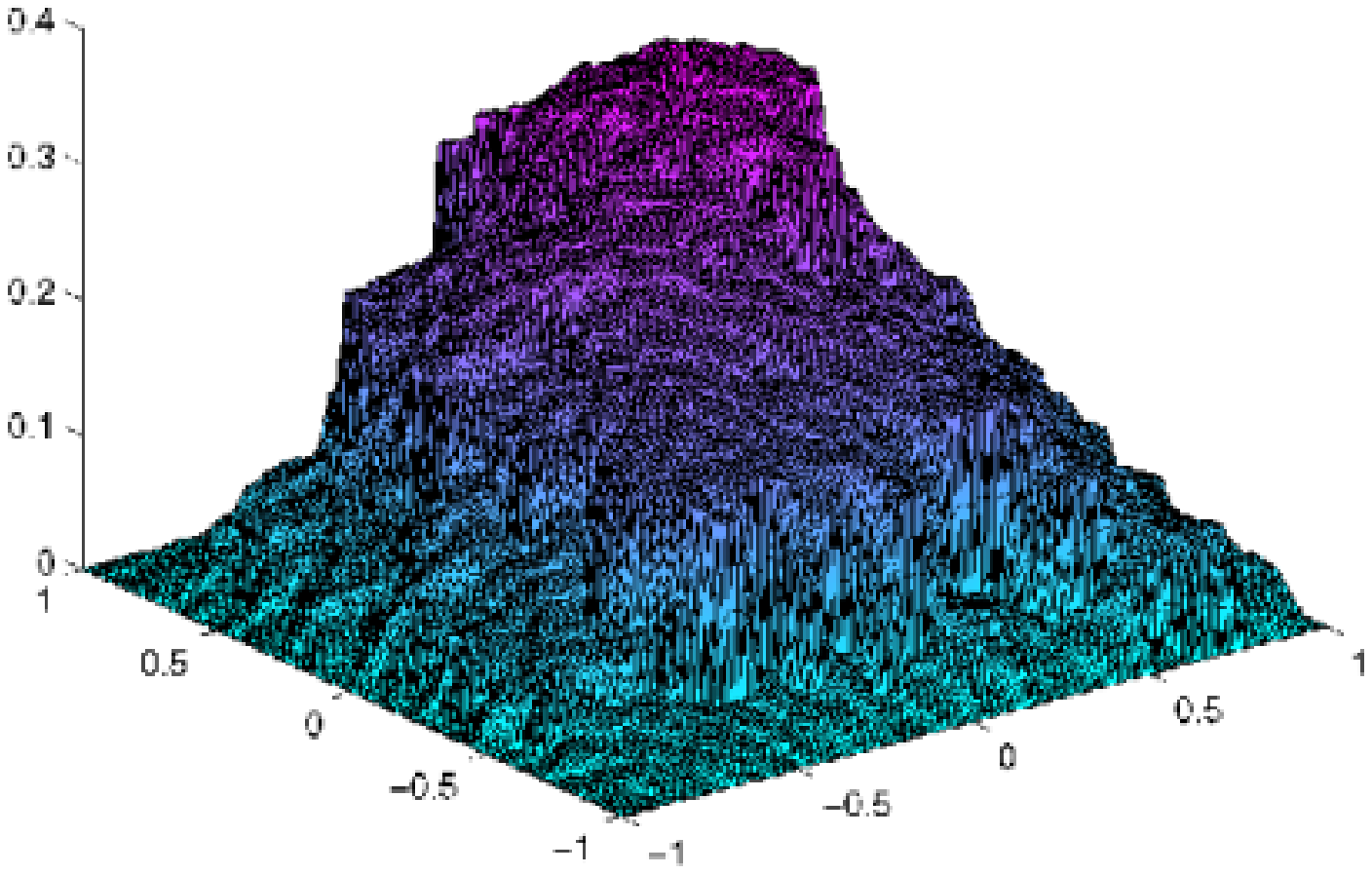}}
    \goodgap
    \subfigure[$v_h$.]
    {\includegraphics[width=0.35\textwidth,height= 0.3\textwidth]{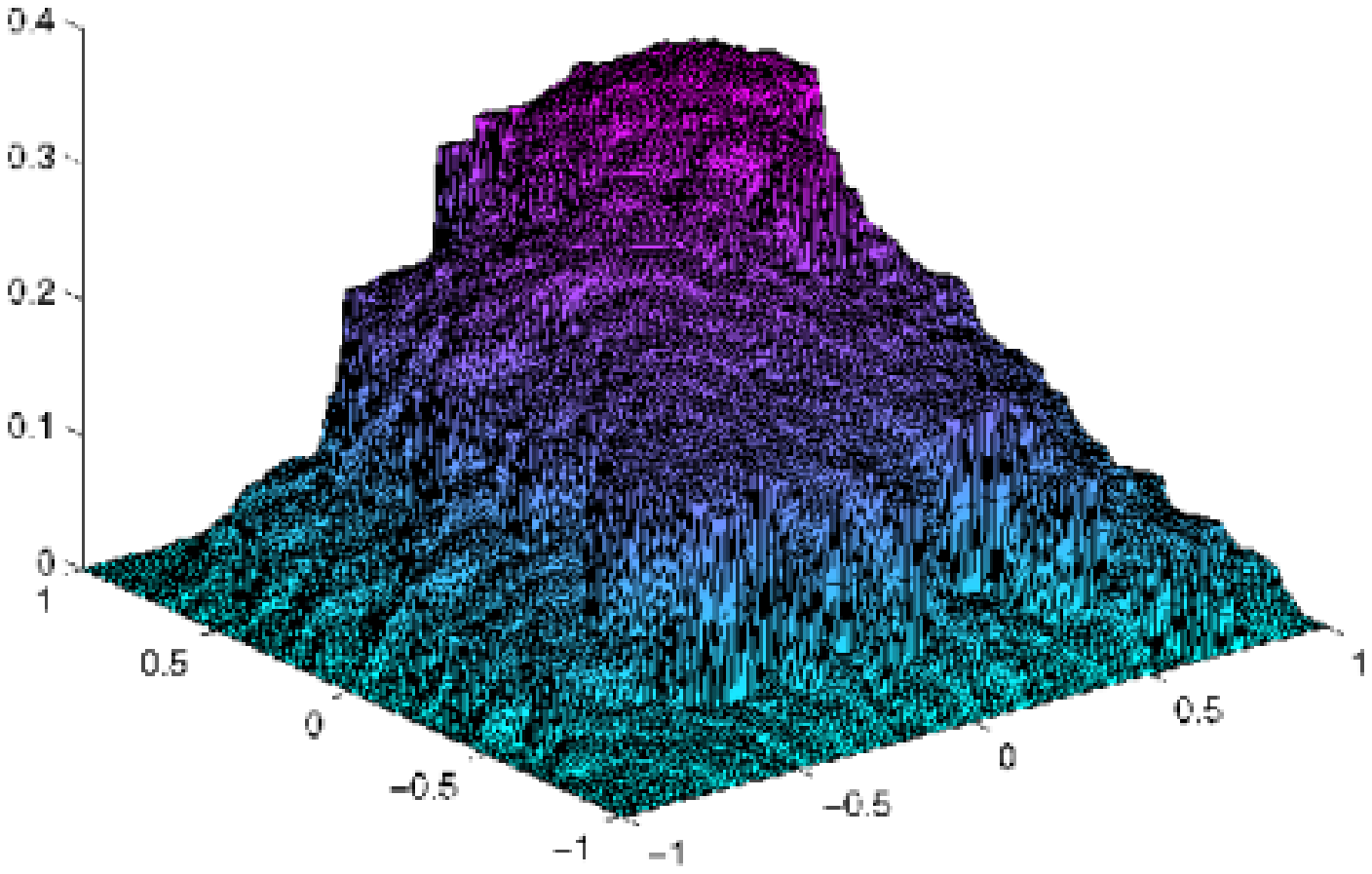}}\\
    \caption{$u$ computed with \emph{dof} $261121$ and $v_h$ computed with
\emph{dof} $9$ at time 1, they are interpolated on a coarser mesh in order to
have a clear picture.}
    \label{eruuhzaa4}
\end{center}
\end{figure}

Figure \ref{eruuhzaa4} shows $u$ computed with $261121$ \emph{dof} and $v_h$ 
computed with $9$ \emph{dof} in the case $g=1$ at time $1$ using method
FEM\_$\psi_{sp}$. They are visually almost the same even in terms of small
scale features. Table \ref{siteperco} gives the numerical errors for $g=1$
with respect to different coarse and fine \emph{dof}.

\begin{table}
\begin{center}
\caption{Example \ref{exa:siteprcoc4}, numerical errors for $g=1$}
\label{siteperco}
\begin{tabular}{|c|c|c|c|c|c|}
\hline $dof_f$&$dof_c$& $L^{1}$& $L^{\infty}$& $L^{2}$& $H^{1}$\\
\hline \multirow{3}{*}{65025} 
& 9& 0.0750& 0.0777& 0.0729& 0.1528 \\
& 49& 0.0301& 0.0351& 0.00298& 0.0779\\
& 225& 0.0135& 0.0147& 0.0133& 0.0333\\
\hline \multirow{3}{*}{261121} 
& 9& 0.0752& 0.0779& 0.0731& 0.1533 \\
& 49& 0.0302& 0.0345& 0.0299& 0.0782\\
& 225& 0.0131& 0.0145& 0.0130& 0.0329\\ \hline
\end{tabular}
\end{center}
\end{table}

\begin{figure}[httb]
\begin{center}
\includegraphics[%
  scale=0.3]{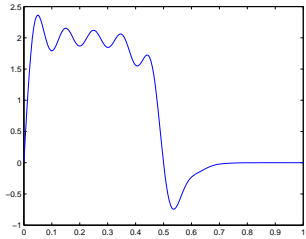}
\caption{$t\rightarrow g(0,t)$} \label{tds}
\end{center}
\end{figure}

Finally, we consider the site percolating medium, with Neumann
boundary condition and a more realistic forcing term. The source 
term is given by $g(x,t)=T(t)X(x,y)$, $X(x,y)$ is the Gaussian source 
function described by
\begin{equation}
X(x,y)=\frac{1}{\sqrt{2\pi\sigma^2}}\exp\big(-\frac{x^2+y^2}{2\sigma^2}
\big),
\end{equation}
with $\sigma=0.05$, $T(t)=T_1(t)T_2(t)$
\begin{equation}
T_1(t)=\sum^{10}_12\frac{1-(-1)^k}{k\pi}\sin(2k\pi t),
\end{equation}
and $T_2(t)=\mbox{erfc}(8(t-0.5))$, $erfc$ is the complementary error
function. We use this source term to emulate a source acting around the origin 
before $t=0.5$, then suddenly decays. See Figure \ref{tds} for $T(t)$ in 
$(0,1)$.

In fact, our future goal is to simulate the response of an explosion, 
usually this is done with a so-called Ricker function \cite{MR1716789} , i.e., 
$g(x,y,t)=\delta_0(x-x_s,y-y_s)R(t)$ with
\begin{equation}
R(t)=(1-2\pi^2(f_0t-1)^2)\exp(-\pi^2(f_0 t-1)^2]
\end{equation}
where $\delta_0$ is the Dirac function and $f_0$ is called the central 
frequency of the source wavelet. It is clear that Ricker function does not
belong to $L^2(\Omega_T)$, our analysis does not apply and numerical experiment
failed in this case. Therefore we would like to test the above modified source 
term first.

Numerical errors for this modified source are given in Table 
\ref{tgsfexa3}. The errors are acceptable but not so good as $g=1$. Possible 
solutions include the adaptive integration in time and adaptation in
space around the source. 

\begin{table}
\begin{center}
\caption{Example \ref{exa:siteprcoc4}, numerical errors for modified source}
\label{tgsfexa3}
\begin{tabular}{|c|c|c|c|c|c|}
\hline $dof_f$&$dof_c$& $L^{1}$& $L^{\infty}$& $L^{2}$& $H^{1}$\\
\hline \multirow{3}{*}{65025} 
& 9& 0.0429& 0.1580& 0.0505& 0.4065 \\
& 49& 0.0261& 0.0965& 0.0324& 0.3030\\
& 225& 0.0185& 0.0879& 0.0226& 0.2718\\
\hline \multirow{3}{*}{261121} 
& 9& 0.0491& 0.1189& 0.0578& 0.4259 \\
& 49& 0.0263& 0.0969& 0.0324& 0.2780\\
& 225& 0.0178& 0.1139& 0.0221& 0.2474\\ \hline
\end{tabular}
\end{center}
\end{table}

\clearpage
\section{Conclusion and Further Remarks}
From above analysis and numerical examples, we observe that good numerical 
approximations can be obtained with much fewer degrees of freedom for 
acoustic wave equation with heterogeneous coefficients, even for the cases
which do not satisfy the Cordes condition. Compared with the multiscale finite element 
method which compute the basis locally, our method has much better accuracy, especially 
for problems with strong non-local effects. 

As it has been done in \cite{OwZh05}, once one understand that the 
key idea for the homogenization of \eref{waveeqn} lies in its higher 
regularity properties with respect to harmonic coordinates  one can homogenize 
\eref{waveeqn} through a different numerical method (such as a finite volume 
method).

Moreover, it could be observed that one could use any set of $n$
linearly independent solutions of \eref{waveeqn} instead of the
harmonic coordinates. The key property allowing the homogenization
of \eref{waveeqn} lies in the fact that if the data (right hand side and 
initial values) has enough integrability then the space of solutions is at 
small scales close in $H^1$ norm to a space of dimension $n$. Thus once one 
has observed at least $n$ linearly independent solutions of \eref{waveeqn}, 
one has seen all of them at small scales. 

Write $L:=-\nabla a\nabla$. $L^{-1}$ maps $H^{-1}(\Omega)$ into
$H_{0}^{1}(\Omega)$, it also maps $L^{2}(\Omega)$ into $V$ a sub-vector
space of $H_{0}^{1}(\Omega)$. The elements of $V$ is close in $H^{1}$
norm to a space of dimension $n$ (the dimension of the physical space
$\Omega$) in the following sense.

Let $\mathcal{T}_{h}$ be a triangulation of $\Omega\subset\mathbb{R}^{n}$
of resolution $h$ (where $0<h<{\tt diam}(\Omega)$). Let $\Lambda$
set of mappings from $\mathcal{T}_{h}$ into the unit sphere of 
$\mathbb{R}^{n+1}$ (if $\lambda\in\Lambda$ then $\lambda$is constant on 
each triangle $K\in\mathcal{T}_{h}$ and $\|\lambda(K)\|=1$), then 
\begin{equation}
\sup_{v_{1},v_{2},\dots,v_{n+1}\in V}\inf_{\lambda\in\Lambda}\frac{\|\sum_{i=1}^{n+1}\lambda_{i}v_{i}\|_{H_{0}^{1}(\Omega)}}{\sum_{i=1}^{n+1}\|\nabla a\nabla v_{i}\|_{L^{2}(\Omega)}}\leq Ch\label{supinf}
\end{equation}

Equation \eref{supinf} is saying that any $n+1$ elements of $V$
are (at an $h$ approximation in $H^{1}$ norm) linearly dependent.
Recall that $n+1$ vectors are linearly dependent in a linear combination
(with non null coefficients) of these vectors in the null vector.
In \eref{supinf} the linear combination of the $n+1$ vectors is at
relative distance of order $h$ (resolution of the triangulation)
from $0$.

We notice that some recent results using global information \cite{AaEfJi07,
JiEfGi07, JiEfGi07b} are formulated in a partition of unity framework 
\cite{BaMe97}. In this case, $\{1,F_1,\cdots,F_n\}$ can be used to construct 
the local approximation space.

\def\cprime{$'$} \def\cprime{$'$}
  \def\polhk#1{\setbox0=\hbox{#1}{\ooalign{\hidewidth
  \lower1.5ex\hbox{`}\hidewidth\crcr\unhbox0}}} \def\cprime{$'$}

\end{document}